\documentclass[reqno,11pt,a4paper]{amsart}
\usepackage{hyperref}
\usepackage{amsmath,amssymb,amsthm,amsfonts}
\usepackage{enumerate}
\usepackage{a4wide}
\usepackage[matrix,arrow,tips]{xy}
\usepackage{h2def}
\newcommand{\restrict}[2][{}]{{\rho_{#1}\left(#2\right)}}

\begin{document}

\date{\today}
\title{Generalised Hilbert Numerators II} 
\author{Jan Snellman}

\address{School of Informatics\\ 
  University of Wales\\
  Dean Street, Bangor \\
  Gwynedd LL57 1UT\\
  Wales\\
  UK
}

\email{jans@matematik.su.se}
\thanks{The author was supported by grants from Svenska Institutet and
  Kungliga Vetenskapsakademin.}

  

\begin{abstract}
    Let \RP\ denote the \emph{large polynomial ring} (in the sense of
  Halter-Koch \cite{HK:genpolalg}) on  the  set
  \(X=\set{x_1,x_2,x_3,\dots}\) of  indeterminates. For each integer \(n\),
  there is a \emph{truncation homomorphism} \(\restrictFunc{n}: \RP \to
  K[x_1,\dots,x_n]\). If 
  \(I\) is a homogeneous ideal of \RP, then the \(\Nat\)-graded Hilbert series 
  of \(\frac{K[x_1,\dots,x_n]}{\restrict[n]{I}}\) can be written as
  \(\frac{g_n(t)}{(1-t)^n}\); it was shown in \cite{Snellman:HilbNumPubl} that 
  if in addition \(I\) is what we call \emph{\lfg}, then \(g_n(t) \to g(t) \in 
  \Z[[t]]\), the so-called \emph{Hilbert numerator} of \(I\).
  
  In this article, we generalise this result to \(\Nat^r\)-graded \lfg\
  ideals. For monomial ideals in \RP, we define the \([X]\)-graded Hilbert
  numerator as the Hilbert numerator of the contracted monomial ideal in
\(K[X]\), for which the standard combinatorial and homological methods
  for calculating multi-graded Hilbert series of monomial ideals in finitely
  many variables apply. Finally, we show that all
  polynomial \(\Nat\)-graded Hilbert numerators can be obtained from ideals
  generated in finitely many variables, and that the the closure in
  \(\Z[[t]]\)  of this set is the set of all
  \(\Nat\)-graded Hilbert numerators.     

  Our main tools are the approximation theorem of \cite{Snellman:TruncPubl},
  relating the initial ideal with the initial ideal of the truncated ideal,
  and the identification of \(K[[X]]\) with the ring of all number-theoretic
  functions \cite{NumThe} which allows the passing from the characteristic
  function of the complement of a monomial ideal to its Hilbert numerator to
  be seen as an example of M\"obius inversion.

\end{abstract}
\maketitle
\tableofcontents
\newcommand{\alfg}{lfg}
\begin{section}{Introduction}
The ring \RP\ of formal polynomials was used by Halter-Koch
\cite{HK:genpolalg} to 
study polynomial functions on modules. The author used it to study initial
ideals of generic ideals of the same type, for instance generated by a
quadratic and a 
cubic generic form, but in ever more variables 
\cite{Snellman:GbRprimPubl,Snellman:TruncPubl,Snellman:Revlexico}. 
In brief, it is the largest
\(\Nat\)-graded subring of the
power series ring \(K[[X]]\) on a countable set \(X\) of indeterminates. There 
are \emph{truncation maps} \(\restrictFunc{n}: \RP \to K[x_1,\dots,x_n]\), and 
inclusion maps going the other way, which relate ideals in \RP\ with 
sequences \((I_n)_{n=1}^\infty\) of ideals, where \(I_n\) is an ideal in
\(K[x_1,\dots,x_n]\), and where \(I_{n+1}\) maps to \(I_n\) under the map 
\begin{displaymath}
  K[x_1,\dots,x_{n+1}] \twoheadrightarrow
  \frac{K[x_1,\dots,x_{n+1}]}{(x_{n+1})} \simeq K[x_1,\dots,x_n]. 
\end{displaymath}
As an example, choose positive integers \(d_1, \dots, d_r\), and let for each
positive integer \(n\)
\(f_1^{(n)},\dots,f_r^{(n)}\) be  forms in \(n\) variables of degree \(d_1\)
to \(d_r\).  Suppose furthermore that \(f_\ell^{(n+1)} - f_\ell^{(n)}\) is
divisible by \(x_{n+1}\), and that the coefficients of the forms are choosen
randomly. Let \(I_n = (f_1^{(n)},\dots,f_r^{(n)})\), and let \(>\) the the
lexicographic term order on the various polynomial  rings. Then the initial
ideals \(\init(I_n)\) will converge to a monomial ideal in infinitely many
indeterminates, as \(n \to \infty\), and this monomial ideal is the
lex-initial ideal of \(I = (f_1,\dots,f_r) \subset \RP\), where \(f_\ell =
\lim f_\ell^{(n)}\) \cite{Snellman:TruncPubl}. In this case, we have that for
\(n \ge r\), the Hilbert 
series of \(\frac{K[x_1,\dots,x_n]}{I_n}\) is \((1-t)^{-n} \prod_{j=1}^r
(1-t^{d_j})\). 

For arbitrary homogeneous, finitely generated ideals \(I \subset \RP\), we
conjecture that \[(1-t)^n H(\frac{K[x_1,\dots,x_n]}{\restrict[n]{I}}; t)\] is
eventually constant, where \(H(\frac{K[x_1,\dots,x_n]}{\restrict[n]{I}}; t)\)
denotes the Hilbert series of
\(\frac{K[x_1,\dots,x_n]}{\restrict[n]{I}}\). What we have in fact proved
\cite{Snellman:HilbNumPubl} is
that for all homogeneous ideal which are countably generated and which have
but finitely many minimal generators of each total degree (we call such an
ideal \emph{\lfg}), the polynomials 
\[(1-t)^n H(\frac{K[x_1,\dots,x_n]}{\restrict[n]{I}}; t) \to p(t) \in
\Z[[t]],\] we call the power series
 \(p(t)\) the \emph{generalised Hilbert numerator}
of \(I\). The outstanding question is thus whether for finitely generated ideals,
\(p(t)\) is always a polynomial.

In this article, we first show that for monomial ideals in \(K[X]\), the
polynomial ring on countably many indeterminates, the usual methods of
calculating multigraded Hilbert series can be used, and that this Hilbert series can
always be written \(\frac{p(X)}{\prod_{i=1}^\infty (1-x_i)}\). We call
\(p(X)\) the \([X]\)-multigraded Hilbert numerator of the ideal. For a
\lfg\ ideal in \(\RP\), we form its initial ideal, contract it to a
monomial ideal in \(K[X]\), calculate the \([X]\)-multigraded Hilbert
numerator of that, then collapse the grading to a \(\Nat\)-grading to get the
\(\Nat\)-graded Hilbert numerator. Apart from providing a more attractive
proof of the existence of Hilbert numerators, this methodology yields
immediately \(\Nat^r\)-graded Hilbert numerators for \(\Nat^r\)-graded \lfg\
ideals of \RP.

We give exact (although opaque) descriptions of the set of all \([X]\)-graded
Hilbert numerators of monomial ideals, and the set of all \(\Nat\)-graded
Hilbert numerators of \lfg\ ideals of \RP. The latter result can be described
briefly as follows: the set of \emph{all polynomial} \(\Nat\)-graded Hilbert
numerators 
is the set of sufficiently high iterated differences of admissible
\(H\)-vectors (in the sense of Macaulays characterisation of admissible
Hilbert functions of finitely generated algebras), and the set of \emph{all}
\(\Nat\)-graded Hilbert numerators is the closure in \(\Z[[t]]\) of the
previous set. 
\end{section}

\begin{section}{Notation}
  Let \(\Z,\Nat,\Nat^+\) denote the set of integers, non-negative integers,
  and positive integers, respectively. For any set \(A\) and any positive
  integer \(k\), we denote by \(\binom{A}{k}\) the set of \(k\)-subsets of
  \(A\), by \([A]\) the free abelian monoid on \(A\), and by \([A]_k\) the
  subset of monomials of total degree \(k\). If \(m \in [A]_k\) we write
  \(\tdeg{m} = k\).

  If \(B\) is another set, then \(B^A\) denotes the set of all functions \(A
  \to B\). If \(A\) is pointed, that is, has an extinguished zero element (for
  instance, \([A]\) is pointed), then \(B^{(A)}\) denotes the set of all
  finitely supported maps \(A \to B\).
  
  If \(K\) is a commutative ring, then \(K^{[A]}\) becomes a commutative
  \(K\)-algebra under component-wise addition and multiplication of scalars,
  and with multiplication given by the \emph{convolution product}
  \begin{equation}
    f \times g (m) = \sum_{\divides{t}{m}} f(t) g(m/t).
  \end{equation}
  We denote this ring by \(K[[A]]\).  The set \(B^{([A])}\) is a subring,
  which we denote by \(K[A]\).
\end{section}

\begin{section}{Rings of formal power series and formal polynomials in
    countably many indeterminates} Let \(K\) be a field of containing the
    rational numbers, and let \(X = \set{x_1,x_2,x_3,\dots}\) be a set of
    indeterminates.  Form the large power series ring \(K[[X]]\) and the
    polynomial ring \(K[X]\) as above. For \(K[[X]] \ni f = \sum_{m \in [X]}
    c_m m\) we define
  \begin{equation}
    \supp(f) = \setsuchas{m}{c_m \neq 0}
  \end{equation}

  If \(t \in [X]\), we define
  \begin{equation}
    [t]f = c_t.
  \end{equation}

  The ring \(K[X]\) is \([X]\)-graded, and in particular \(\Nat\)-graded,
  whereas \(K[[X]]\) is not. The largest \([X]\)-graded subring of \(K[[X]]\)
  is \(K[X]\), whereas the largest \(\Nat\)-graded subring is the ring \RP\
  generated by all \emph{bounded elements}: an element \(f \in K[[X]]\) is
  bounded if
  \[\tdeg{f} := \sup \left( \setsuchas{\tdeg{m}}{m \in \supp(f)}\right) < \infty.\]
  Another way of putting this is the following.
    \begin{definition}
    We define the \emph{total-degree} filtration on \(K[[X]]\) and its various 
    subrings by 
    \begin{equation}
      \label{eq:tdegfit}
      \tdegfilt{d} K[[X]] = \setsuchas{f \in K[[X]]}{\tdeg{f} \le d}
    \end{equation}
    For \(K[[X]] \ni f = \sum_{m \in [X]} c_m m\), we put 
    \begin{equation}
      \label{eq:tdf}
      \tdegfilt{d} f = \sum_{\substack{m \in [X] \\ \tdeg{m} \le d}} c_m m
    \end{equation}
  \end{definition}
  Then \(\RP = \cup_{d \ge 0} \tdegfilt{d} K[[X]]\).
  
  It is shown in \cite{Snellman:GbRprimPubl} that \RP\ is also the maximal
  subring of \(K[[X]]\) with the following property: given any multiplicative
  total order \(>\) on \([X]\) whose restriction to \([X]_1\) is
  order-isomorphic to \(-\omega\) (such a \(>\) will be called a \emph{term
  order} on \([X]\)), the support of any non-constant element \(f\) contains a
  maximal element \(\init_>(f)\). Putting \(\init_>(1)=0\), \(\init_>(0)=
  -\infty\), we can regard
  \begin{displaymath}
    \init_>: \RP \to [X] \cup \set{-\infty}
  \end{displaymath}
  as a \([X]\)-valuation, which induces a \([X]\)-filtration of \RP\ by
  \begin{equation}
    \begin{split}
      \mathfrak{F}^{<m} \RP &= \setsuchas{f \in \RP}{\init_>(f) < m} \\
      \mathfrak{F}^{\le m} \RP &= \setsuchas{f \in \RP}{\init_>(f) \le m}
    \end{split}  
  \end{equation}
  We then have a canonical map
  \begin{equation}
    \begin{split}
      \RP & \to \gr(\RP) = \bigoplus_{m \in [X]} \frac{\mathfrak{F}^{\le m}
      \RP}{\mathfrak{F}^{< m} \RP} \simeq K[X] \\ f & \longmapsto \init_>(f)
    \end{split}  
  \end{equation}
  This map sends an ideal \(I \subset \RP\) to its \emph{initial ideal}
  \begin{equation}
    \init_>(I) = \RP \setsuchas{\init_>(f)}{f \in I}
  \end{equation}
  which is a \emph{monomial ideal}, that is, generated by monomials. We note
  that
  \begin{enumerate}
  \item Every monomial ideal is its own initial ideal,
  \item Extension and contraction of ideals gives a bijection between monomial
ideals in \(K[X]\), \RP\ and \(K[[X]]\),
  \item Monomial ideals in \(K[X]\), \RP\ or \(K[[X]]\) correspond bijectively
to monoid ideals in \([X]\).
  \end{enumerate}
  Because of this identification, we shall say that a monoid ideal has a
  certain property whenever the corresponding monomial ideal has.

  \begin{theorem}[Snellman \cite{Snellman:GbRprimPubl}]\label{thm:lfg}
    For a \(\Nat\)-graded ideal \(I \subset \RP\), the following are
    equivalent:
    \begin{enumerate}
    \item \(I\) is generated by a \emph{locally finite} set, that is a set
      \begin{equation}
        F = \bigcup_{d=1}^\infty F_i, \quad \forall i: F_i \in \RP_i, \quad
        \forall i: \# F_i < \infty
      \end{equation}
    \item For each positive integer \(d\),
      \begin{equation}
        \dim_K \left( \frac{I_d}{\sum_{i=1}^{d-1} \RP_i I_{d-i}} \right) <
        \infty
      \end{equation}
    \end{enumerate}
  \end{theorem}
  We call such ideals \emph{\lfg} (\alfg). By our previous remark, we can talk about
  \alfg\ monoid ideals, as well.

  \begin{theorem}[Snellman \cite{Snellman:GbRprimPubl}]
    Let \(>\) be a term-order on \([X]\), and \(I \subset \RP\) a homogeneous
    ideal. Then \(I\) is \alfg\ if and only if \(\init_>(I)\) is.
  \end{theorem}

  \begin{subsection}{Inverse limits}

    We shall need to relate elements in \RP\ with their \emph{truncations} in
    \(K[X_n]\). The necessary machinery is as follows.

    For any positive integer \(n\), we put \(X_n = \set{x_1,\dots,x_n}\), and
    let \([X_n]\) be the free abelian monoid on \(X_n\). We define the
    polynomial ring \(K[X_n]\) and the power series ring \(K[[X_n]]\) as
    above. For \(i < j \) there is a commutative diagram of \(K\)-multilinear
    maps
    \begin{equation}
      \xymatrix{ K[X] \ar [rrr] \ar [dr] \ar [dd] &&& K[[X]] \ar [dl] \ar
      [dd]\\ & K[X_j] \ar [r] \ar [dl]& K[[X_j]] \ar [dr] \\ K[X_i] \ar
      [rrr]&&& K[[X_i]] }
    \end{equation}
    with the horisontal arrows given by inclusions, and the remaining ones
    given by (restrictions of) the truncation maps
    \begin{equation}
      \begin{split}
        \restrictFunc{n}: [X] &\to [X_n] \cup \set{0} \\ 
        m & \mapsto
        \begin{cases}
          m & m \in [X_n] \\ 0 & m \not \in [X_n]
        \end{cases}
        \\
        \restrictFunc{n}: K[[X]] &\to K[[X_n]] \\ 
        \sum_{m \in [X]} c_m m & 
        \mapsto \sum_{m \in [X]} c_m \restrict[n]{m},
      \end{split}
    \end{equation}

    With respect to these inverse systems, we have that \(\varprojlim K[[X_n]]
    \simeq K[[X]]\), whereas
    \begin{displaymath}
      K[X] \subsetneq \RP \subsetneq \varprojlim K[X_n] \subsetneq K[[X]].
    \end{displaymath}
    In fact,
    \begin{equation}
      \begin{split}
        \varprojlim K[X_n] & = \setsuchas{f \in K[[X]]}{\forall n: \,
        \restrict[n]{f} \in K[X_n]} \\ \RP & = \setsuchas{f \in\varprojlim
        K[X_n]}{f \text{ is bounded }}.
      \end{split}  
    \end{equation}
    Furthermore, the ring \(\varprojlim K[X_n]\) is endowed with a natural
    topology, the inverse limit topology (where all \(K[X_n]\) are discrete),
    and the ring \RP\ is a dense subring. The topology on \(\RP\) can be
    characterised by giving the closure of an arbitrary subset \(A \subset
    \RP\):
    \begin{equation}
      \bar{A} = \setsuchas{f \in \RP}{\forall n: \, \restrict[n]{f} \in
      \restrict[n]{A}}.
    \end{equation}
    It was proved in \cite{Snellman:Ringt1Publ} that with respect to this
    topology, \alfg\ ideals in \RP\ are closed. It was also proved that the
    closed monomial ideals are precisely the \alfg\ monomial ideals.

  \end{subsection}

  \begin{subsection}{Topologies on the set of ideals of \RP, and a
      ``continuity'' result}
    
    \begin{definition}
      Let \(\ideal{}, \ideal{c}, \ideal{h}, \ideal{l}, \ideal{m}\) denote the
      following sets of ideals in \RP: all ideals, closed ideals, homogeneous
      ideals, \alfg\ ideals, monomial ideals. We will also use combinations of
      letters to denote intersections, for instance
      \begin{displaymath}
        \ideal{lm} = \ideal{l} \cap \ideal{m}
      \end{displaymath}
      denotes the \alfg\ monomial ideals.
    \end{definition}

    \begin{proposition}\label{prop:top}
      \begin{enumerate}[(i)]
      \item The function
        \begin{equation}\label{eqn:metricdd}
          d(I,J) = 2^{-n}, \qquad n = \max \setsuchas{n}{\restrict[n]{I} =
          \restrict[n]{J}}
        \end{equation}
        gives a metric on \ideal{c}.
      \item The function
        \begin{equation}\label{eqn:metricdn}
          \hat{d}(I,J) = 2^{-d}, \qquad d = \max \setsuchas{d}{\tdegfilt{d} I
          = \tdegfilt{d}J}
        \end{equation}
        gives a metric on \ideal{h}.
      \item Define a convergence structure on \ideal{m} by dictating that
        \(I_n \tostar I \in \ideal{m}\) if and only if,
        \begin{equation}\label{eqn:convweak}
          \forall m \in [X]:
          \exists N(m) \in \Nat^+: \forall n>N(m): \, m \in I \iff m \in I_n
        \end{equation}
        Then the corresponding topology is weaker than both the previous
        topologies.
      \end{enumerate}
    \end{proposition}
    \begin{proof}
      \begin{enumerate}[(i)]
      \item It is clear the \(d(I,J) = d(J,I) \ge 0\). Since \(I,J\) are
closed, \(d(I,J)=0\) if and only if \(I=J\). If \(A,B,C\) are closed ideals,
and if \(d(A,B) \le 2^{-n}\), \(d(B,C) \le 2^{-n}\), then \(\restrict[n]{A} =
\restrict[n]{B} = \restrict[n]{C}\), hence \(d(A,C) \le 2^{-n}\). Thus the
triangle inequality holds.
      \item Obvious.
      \item Let \(m \in [X_v]_d\), let \(I,I_1,I_2,I_3,\dots\) be monomial
ideals, and suppose that either \(d(I_n,I) \to 0\) or \(\hat{d}(I_n,I) \to
0\). In the first case, there is an \(N(v)\) such that \(\restrict[v]{I_n} =
\restrict[v]{I}\) whenever \(n \ge N(v)\): since \(m \in I \iff m \in
\restrict[v]{I}\), and similarly for \(I_n\), it follows that \(m \in I\) if
and only if \(m \in I_n\) for \(n \ge N(v)\).  In the second case, there is an
\(\hat{N}(d)\) such that \(\tdegfilt{d}I_n = \tdegfilt{d}I\) whenever \(n \ge
\hat{N}(d)\): since \(m \in I \iff m \in \tdegfilt{d}I\), and similarly for
\(I_n\), it follows that \(m \in I\) if and only if \(m \in I_n\) for \(n \ge
\hat{N}(d)\).
      \end{enumerate}
    \end{proof}

    \begin{theorem}
      Let \(>\) be a term-order on \([X]\). Then the map
      \begin{equation}\label{eqn:inmap}
        \begin{split}
          \init_>: \ideal{l} & \to \ideal{lm} \\ I & \mapsto \init_>(I)
        \end{split}
      \end{equation}
      is continuous with respect to the \(\hat{d}\)-metric.  If \(>\) is the
      degree-reverse lexicographic term order, then 
      \begin{equation}
        \label{eq:swap}
        \forall n: \qquad
        \restrict[n]{\init_>(I)} = \init_>(\restrict[n]{I})
      \end{equation}
      from which it follows that \eqref{eqn:inmap} is continuous with respect
      to the \(d\)-metric.
    \end{theorem}
    \begin{proof}
      Using the results of \cite{Snellman:GbRprimPubl}, it is straight-forward
      to show that if \(I,J\) are \alfg\ ideals such that \(\tdegfilt{d} I =
      \tdegfilt{d} J\), then for any term-order \(>\), \(\tdegfilt{d}
      \init_>(I) = \tdegfilt{d} \init_>(J)\). Hence the first result follows.
      
      It is immediate that the identity \eqref{eq:swap} implies continuity of
      \eqref{eqn:inmap}. For all term orders, the LHS of \eqref{eq:swap} is
      included in the RHS, so we need to prove that the reverse inclusion
      holds for the degree-reverse lexicographic term order. Let \(f \in I\)
      be homogeneous of degree \(d\); then the monomials in \(\supp(f)\) are
      ordered as follows: first the ones in \([X_1]_d \cap \supp(f)\), if any,
      then the ones in \(([X_2]_d \setminus [X_1]_d) \cap \supp(f)\), and so
      on. Let \(m = \init(f)\), then \(\init(\restrict[n]{f}) = m\) for \(n\)
      sufficiently large, and \(0\) otherwise. In the same way,
      \(\restrict[n]{m} = m\) for sufficiently large \(n\), and \(0\)
      otherwise. So \(\mathrm{RHS} \ni \init(\restrict[n]{f}) = \restrict[n]{\init(f)}
      \in \mathrm{LHS}\).
    \end{proof}
    
    The following result is a key one: it is what will allow us to define
    Hilbert numerators of \alfg\ ideals by passing to their initial ideals.
    \begin{theorem}[Snellman \cite{Snellman:TruncPubl}]  \label{thm:inconv}
      If \(I\) is a \alfg\ ideal, and \(>\) is a term-order on \([X]\), then
      \begin{equation}
        \label{eq:inreco}
        \hat{d}(\init_>(\restrict[n]{I})^e,\init_>(I)) \to 0 \quad \text{ as
          } n \to \infty.
      \end{equation}
    \end{theorem}
  \end{subsection}
\end{section}

\begin{section}{Monoid ideals and arithmetic on
    \protect\(\protect\Z[[X]]\protect\)} 
  \begin{subsection}{Topologies on \protect\(\protect\Z[[X]]\protect\)}
    Unless otherwise stated, we henceforth assume that \(\Z[[X]] = \Z^{[X]}\)
    is given the product topology. With this topology, \(f_n \to f\) if for
    all \(m \in [X]\), there is an \(N(m)\) so that for \(n \ge N(m)\), \([m]f
    = m[f_n]\). An infinite sum \(\sum_{n} f_n\) is convergent if and only if
    each monomial \(m \in [X]\) occurs in but finitely many of the sets
    \(\supp(f_n)\).

    We can also topologise \(\Z[[X]]\) by means of the total degree
    filtration: a sequence \((f_n)_{n=1}^\infty\), \(f_n \to f\) if and only
    if
    \begin{displaymath}
      \forall d \in \Nat: \, \exists N(d) \in \Nat^+: \, \forall v > N(d):
      \forall m \in [X]_d: [m]f_v = [m]f.
    \end{displaymath}
    This is a stronger topology than the previous one. We shall use it in
    particular for the study of the subring \(\mathcal{S}\), to be defined
    later. For later use, we note the following:
    \begin{lemma}\label{lemma:transclos}
      The total degree filtration topology on \(\Z[[X]]\) gives a linear
      topology, and hence additiv translation with arbitrary elements, and
      multiplicative translation with invertible elements, are closed mappings.
    \end{lemma}
    \begin{proof}
      We put \[J_d = \set{0} \cup \setsuchas{f \in \Z[[X]]}{\supp(f) \subset
        \cup_{v=d} [X]_v}.\]
      Then the \(J_d\)'s are clopen ideals which form a fundamental system of
      neighbourhoods of zero.

      It follows \cite{Bourbaki:Commutative} that for any subset \(A \subset
      \Z[[X]]\), the closure is given by 
      \begin{displaymath}
        \bar{A} = \cap_{d=1}^\infty (A + J_d).
      \end{displaymath}
      Hence if \(h \in \Z[[X]]\) and \(A \subset \Z[[X]]\), 
      then
      \begin{displaymath}
        \bar{h+A} = \cap_{d=1}^\infty (h+A + J_d) = 
        h + \cap_{d=1}^\infty (A  + J_d)  = h +\bar{A},
      \end{displaymath}
      where the crucial inclusion 
      \begin{math}
        \cap_{d=1}^\infty (h+A + J_d)  \subset 
        h + \cap_{d=1}^\infty (A  + J_d)
      \end{math}
      is proved as follows. If \(f \in  h +A + J_d\) for all \(d\), 
     then \(f - h \in A + J_d\) for all \(d\), hence \(f -h \in
     \cap_{d=1}^\infty A + J_d\), hence \(f \in h + \cap_{d=1}^\infty A +
     J_d\). 

     If \(h\) has a multiplicative inverse \(h^{-1}\), then 
     \begin{displaymath}
        \bar{hA} = \cap_{d=1}^\infty (hA + J_d) = 
        h  \cap_{d=1}^\infty (A  + J_d)  = h \bar{A};
      \end{displaymath}
      the inclusion 
      \begin{math}
        \cap_{d=1}^\infty (hA + J_d) \subset 
        h  \cap_{d=1}^\infty (A  + J_d)
      \end{math}
      is proved as follows. Suppose that \(f \in  h A + J_d\) for all \(d\), 
      then \(h^{-1}f  \in A + h^{-1}J_d \subset A + J_d\) for all \(d\),
      hence 
      \(h^{-1}f  \in \cap_{d=1}^\infty (A  + J_d)\), hence
      \(f \in h \cap_{d=1}^\infty (A  + J_d)\).
    \end{proof}
  \end{subsection}

  \begin{subsection}{The ring of number-theorethic functions}
    Define \(\Gamma\) to be the set of all maps \(\Nat^+ \to \Z\). With
    component-wise addition and multiplication by scalars, and with the
    \emph{Dirichlet convolution}
    \begin{equation}
      f * g (n) = \sum_{\divides{k}{n}} f(k) g(n/k),
    \end{equation}
    \(\Gamma\) becomes a commutative ring, often referred to as \emph{the ring
    of number-theoretic functions} \cite{NumThe}. The well-known isomorphism,
    given by unique factorisation of integers, between the multiplicative
    monoid \((\Nat^+, \cdot)\) of the positive integers and a denumerable sum
    of copies of \((\Nat,+)\), induces an isomorphism
    \begin{equation}
      \label{eq:theiso}
      \begin{split}
        \Gamma & \to \Z[[X]] \\ f & \mapsto \sum_{m = x_1^{\alpha_1} \cdots
        x_n^{\alpha_n} \in [X]} f(p_1^{\alpha_1} \cdots p_n^{\alpha_n}) m,
      \end{split}    
    \end{equation}

    Define the elements \(\nu, \mu\ \in \Z[[X]]\) by
    \begin{equation}\label{eqn:numu}
      \begin{split}
        \nu &= \sum_{m \in [X]} m  \\
        \mu &=  \prod_{i=1}^\infty (1-x_i) = 
        1 - \sum_{i=1}^\infty
        x_i + \sum_{i < j} x_i x_j - \sum_{i < j < k} x_i x_j x_k + \cdots.
      \end{split}
    \end{equation}
    Then the image of \(\mu\) in \(\Gamma\) is the well-known M\"obius
    function, and M\"obius inversion can be expressed by the formula
    \begin{equation}
      \nu \mu = 1.
    \end{equation}

    We note that we can write
    \begin{equation}
      \begin{split}
        \nu = 1+\sum_{i=1}^\infty \nu_i, & \qquad \mu = 1+\sum_{i=1}^\infty
        \mu_i \\ \nu_i = \sum_{m \in [X]_i} m, & \qquad \mu_i = (-1)^{i}
        \sum_{\sigma \in \binom{[X]}{i}} \sigma
      \end{split}
    \end{equation}
    where \(\nu_i\) is the \(i\)'th \emph{complete symmetric function}
    \cite{SymFunc} and \((-1)^i \mu_i\) is the \(i\)'th \emph{elementary
    symmetric function}.
  \end{subsection}

  \begin{subsection}{Characteristic/generating functions of monoid ideals}
    \begin{subsubsection}{Definitions}
      \begin{definition}
        If \(I\) is a monoid ideal in \([X]\) then
        \begin{equation}
          W(I) = I \setminus \left([X] \setminus \set{1} \right) I
        \end{equation}
        denote the canonical set of minimal generators of \(I\).  We define
        \begin{align}
          \Characteristic(I) &= \sum_{m \in I} m \\ 
          w(I) &= \sum_{m \in W(I)} m \\
          q(I) &= \nu - \Characteristic(I) \\ p(I) &= \mu q(I)
        \end{align}
        We call \(\Characteristic(I)\) the \emph{characteristic function} of
        \(I\), \(q(I)\) the \([X]\)-graded \emph{Hilbert series} of \(I\), and
        \(p(I)\) the \([X]\)-graded \emph{Hilbert numerator} of \(I\).

        For a monomial ideal \(J\) in \(K[X]\) or \RP, we put
        \begin{equation}
          \begin{split}
          \Characteristic(J) &= \Characteristic(J \cap [X])\\
          \quad w(J) &= w(J  \cap [X]) \\
          \quad q(J) &= q(J \cap [X]) \\
          p(J) &= p(J \cap [X]).         
          \end{split}
        \end{equation}

        Similarly, if \(n\) is a positive integer, and if \(I\) is a monoid
        ideal in \([X_n]\), then we put
        \begin{align*}
          \Characteristic^n(I) &= \sum_{m \in I} m \\ q^n(I) &= \sum_{m \in
          [X_n] \setminus I} m \\ p^n(I) &= \restrict[n]{\mu} q^n(I).
        \end{align*}
      \end{definition}
      
      \begin{remark}
      \(\Characteristic(I)\)
      and \(q(I)\) are the \([X]\)-graded Hilbert series of \(I\), regarded as 
      a monomial ideal in \(K[X]\), and
      \(\frac{K[X]}{I}\), respectively. 
      However, the ring \RP\ is not \([X]\)-graded, so in order to attach a
      meaning to \(q(I)\) for a monomial ideal we regard it as a
      limit of the Hilbert  
      series of \(\restrict[n]{I}\), that is, as a limit of
      \(q^n(I)\). 
      \end{remark}

      We note that \(\Characteristic(I)\), \(w(I)\), \(q(I)\), and \(p(I)\)
      all lie in \(\Z[[X]]\).
    \end{subsubsection}

    \begin{subsubsection}{Distributiveness properties}
      \begin{proposition}\label{prop:chidist}
        Suppose that \(I,I_1,I_2,I_3, \dots\) are monomial ideals, and suppose
        that
        \begin{equation}
          \label{eq:condsum}
          \sum_{i=1}^\infty I_n = I,
        \end{equation}
        and that the sum is convergent with respect to the \(\tostar\)
        topology. 
        Then
        \begin{equation}
          \label{eq:chidist}
          \Characteristic(I) = \sum_{i} \Characteristic(I_i) - \sum_{i < j}
          \Characteristic(I_i \cap I_j) + \sum_{i < j < k} \Characteristic(I_i
          \cap I_j \cap I_k) - \cdots,
        \end{equation}
        and the sum is convergent (with respect to the product topology on
        \(\Z[[X]]\)).
        
        Putting \(\hat{p}(I) = p(I) - 1\), we also have that
        \begin{equation}
          \label{eq:pdist}
          \hat{p}(I) = \sum_{i} \hat{p}(I_i) - \sum_{i < j} \hat{p}(I_i \cap
          I_j) + \sum_{i < j < k} \hat{p}(I_i \cap I_j \cap I_k) - \cdots,
        \end{equation}
        and this is a convergent sum.
      \end{proposition}
      \begin{proof}
        If we identify monomial ideals with their characteristic functions,
        and write \(\wedge\) for intersections of ideals, and \(\vee\) for sum
        of ideals, then \(\wedge\) and \(\vee\) correspond to component-wise
        minimum and maximum, and \eqref{eq:chidist} to the identity
        \begin{equation}
          \label{eq:funcdist}
          \bigvee_{i =1}^\infty f_i = \sum_{i} f_i - \sum_{i < j} f_i \wedge
          f_j + \cdots,
        \end{equation}
        where the sum is component-wise. The LHS of \eqref{eq:funcdist} is
        always defined; for the RHS to be defined, it is necessary and
        sufficient that
        \begin{displaymath}
          \forall m \in [X]: \exists N(M): \forall n > N(M): \quad f_i(m) = 0.
        \end{displaymath}
        If this holds, then denoting by \(S\) the cardinality of the finite
        subset \(\setsuchas{j \in \Nat^+}{f_j(m) \neq 0}\), the formula
        \eqref{eq:funcdist} becomes
        \begin{displaymath}
           S - \binom{S}{2} + \binom{S}{3} - \cdots = 
           \begin{cases}
             0 & S = \emptyset \\
             1 & S \neq \emptyset
           \end{cases}
        \end{displaymath}
        a well-know binomial identity.

        To prove \eqref{eq:pdist}, note that \(\hat{p}(I_i) = - \mu
        \Characteristic(I_i)\), hence from \eqref{eq:funcdist} we get that
        \begin{align*}
          \hat{p}(I) &= \mu \Characteristic(I) \\
          &= -\mu \left( \sum_{i}
          \Characteristic(I_i) - \sum_{i < j} \Characteristic(I_i \cap I_j) +
          \sum_{i < j < k} \Characteristic(I_i \cap I_j \cap I_k) - \cdots
          \right) \\ &= \sum_{i} (-\mu \Characteristic(I_i)) - \sum_{i < j}
          (-\mu \Characteristic(I_i \cap I_j)) + \cdots \\ & = \sum_{i}
          \hat{p}(I_i) - \sum_{i < j} \hat{p}(I_i \cap I_j) + \sum_{i < j < k}
          \hat{p}(I_i \cap I_j \cap I_k) - \cdots
        \end{align*}
      \end{proof}
    \end{subsubsection}

    \begin{subsubsection}{Inclusion-exclusion for Hilbert numerators}
      \begin{theorem}\label{thm:inex}
        Let \(I \subset [X]\) be a monoid ideal. If \(\sigma \subset W(I)\) is
        finite, let \(\lcm(\sigma)\) be the least common multiple of the
        elements in \(\sigma\), and let \(\# \sigma\) be the cardinality of
        \(\sigma\). Then
        \begin{equation}\label{eqn:sigma}
          p(I) = \sum_{\sigma} (-1)^{(\# \sigma)} \lcm(\sigma),
        \end{equation}
        where the sum is over all finite subsets of \(W(I)\). Alternatively,
        \begin{equation}\label{eqn:sigmasum}
          p(I) = 1 - \sum_{m \in W(I)} m + \sum_{\sigma \in \binom{W(I)}{2}}
          \lcm(\sigma) - \sum_{\sigma \in \binom{W(I)}{3}} \lcm(\sigma) +
          \cdots
        \end{equation}
      \end{theorem}
      \begin{proof}
        We have that \(I = \sum_{m \in W(I)} (m)\), and that \((m_i) \cap
        (m_j) = (\lcm(m_i,m_j))\).  Hence the result follows from
        \eqref{eq:pdist}, once we have proved that that \(p((m)) = 1 - m\) for
        all \(m \in [X]\). But \(\Characteristic((m)) = \sum_{\divides{m}{t}}
        t\), hence by M\"obius inversion
        \begin{displaymath}
          p((m)) = 1 - \mu \Characteristic((m)) = 1 - \mu
          \left(\sum_{\divides{m}{t}} t\right) = 1-\sum_{\divides{m}{t}} \mu t
          = 1-m.
        \end{displaymath}
      \end{proof}
    \end{subsubsection}

    \begin{subsubsection}{Homology methods}
      \begin{lemma}\label{lemma:approx}
        Let \(I \subset [X]\) be a monoid ideal. Then
        \begin{equation}
          \label{eq:app}
          \forall n \in \Nat^+: \qquad \restrict[n]{p(I)} =
          p^n(\restrict[n]{I})
        \end{equation}
      \end{lemma}
      \begin{proof}
        \begin{align*}
          p^n(\restrict[n]{I}) &= \restrict[n]{\mu} q^n(\restrict[n]{I}) \\ &=
          \restrict[n]{\mu} \sum_{m \in [X_n] \setminus \restrict[n]{I}} m \\
          &= \restrict[n]{\mu} \restrict[n]{\sum_{m \in [X] \setminus I} m} \\
          & = \restrict[n]{\mu} \restrict[n]{q(I)} = \restrict[n]{\mu q(I)} =
          \restrict[n]{p(I)}.
        \end{align*}
      \end{proof}
      Using this lemma, we can immediately extend the various homological
      methods for getting the multigraded Hilbert series of monoid ideals in
      \([X_n]\) (see \cite{MID,LCMlattice}) to work for monoid ideals in
      \([X]\).

      We get
      \begin{theorem}\label{thm:euchar}
        Let \(I \subset [X]\) be a monoid ideal, and let \(m \in [X]\).  Let
        \(\Delta_m =\Delta_m(I) \subset 2^{(\Nat^+)}\) be the following
        simplicial complex:
        \begin{equation}\label{eq:sc}
          \sigma = \set{\sigma_1,\dots,\sigma_r} \in \Delta_m \qquad \iff
          \qquad \frac{m}{\prod_{i=1}^r x_{\sigma_i}} \in I.
        \end{equation}
        Then \(\Delta_m\) is finite, and
        \begin{equation}
          \label{eq:coeffm}
          [m]p(I) =  \tilde{\chi}(\Delta_m(I)) =  \sum_{F \in \Delta_m(I)}
          (-1)^{\tdeg{F} - 1} =  \sum_{i=-1}^{\infty} (-1)^i \dim
          H_i(\Delta_m, K),
        \end{equation}
        where \(\tilde{\chi}\) denotes the the reduced Euler characteristic of
        an abstract 
        simplicial complex, counting the empty set as a \(-1\)-face.
      \end{theorem}

      \begin{theorem}\label{thm:lcmlattice}
        Let \(I \subset [X]\) be a monoid ideal, let \(W=W(I)\) be its minimal
        set of generators, and let \(L_I\) be the lattice of all finite lcm's
        of elements in \(W\), ordered by divisibility. Let \(\hat{0}\) denote
        the minimal element in \(L_I\), and let, for \(m \in L_I\),
        \(\mu(\hat{0}, m)\) denote the value of the M\"obius function of the
        poset \(L_I\), evaluated on the interval \([\hat{0},m]\). 
        Let \(\Delta(\hat{0},m)\) denote the abstract simplicial complex of all
        chains in \((\hat{0},m)\). Then
        we have:
        \begin{equation}
          \label{eq:lcmlat}
          \forall m \in [X]: \qquad [m]p(I) 
          \quad = \quad
          \begin{cases}
            0 & m \not \in L_I \\
            \tilde{\chi}(\Delta(\hat{0}, m)) = 
            \mu(\hat{0},m) & m \in L_I
          \end{cases}
        \end{equation}
      \end{theorem}
      \begin{proof}
        It follows from \cite{LCMlattice} that \(c_m = 0\) for \(m \not \in
        L_I\), and that \(c_m = \tilde{\chi}(\Delta(0,m))\) for \(m \in
        L_I\). By \cite{Stanley:En1}, 
        \begin{math}
          \tilde{\chi}(\Delta(0,m)) = \mu(\hat{0},m)
        \end{math}
        whenever \(m \in L_I\).
      \end{proof}
    \end{subsubsection}

    \begin{subsubsection}{Classifications}
      \begin{proposition}\label{prop:pcond}
        Let \(f =\sum_{m \in [X]} c_m m \in \Z[[X]]\). Then \(f \in
        p(\ideal{m})\) 
        if and only if the following conditions hold:
        \begin{enumerate}
        \item \(\forall m \in [X]: \quad \sum_{\divides{s}{m}} c_s \in
          \set{0,1}\),
        \item If \(\sum_{\divides{s}{m}} c_s =1\) and \(\divides{t}{m}\) then
          \(\sum_{\divides{s}{t}} c_s = 1\).
        \end{enumerate}
      \end{proposition}
      \begin{proof}
        Suppose that \(I\) is a monoid ideal in \([X]\), then \(q(I) = \nu -
        \Characteristic(I)\) is the characteristic function of \(I^c = [X] \setminus 
        I\). This is an \emph{order ideal}, that is, if \(m \in I^c\) and
        \(\divides{t}{m}\), then \(t \in I^c\). It follows that the set of
        \(q(I)\)'s is the set of \(g =\sum_{m \in [X]} d_m m \in \Z[[X]]\)
        such that 
        \begin{enumerate}
        \item \(\forall m \in [X]:  d_m \in  \set{0,1}\),
        \item If \(d_m = 1\) and \(\divides{t}{m}\) then \(d_t = 1\).
        \end{enumerate}
        Since \(p(I) = \mu q(I)\), the result follows by M\"obius inversion.
      \end{proof}

      \begin{proposition}\label{prop:varmax}
        Let \(f =\sum_{m \in [X]} c_m m \in \Z[[X]]\) be the \([X]\)-graded
        Hilbert numerator of a monoid ideal. Let \(m = x_1^{\alpha_1} \cdots
        x_n^{\alpha_n}\). Then
        \begin{equation}
          \label{eq:euc}
          \mathrm{abs}(c_m) \le \binom{n-1}{\lfloor \frac{n-1}{2} \rfloor}
        \end{equation}
      \end{proposition}
      \begin{proof}
        From Theorem~\ref{thm:euchar} we have that \(c_m\) is the reduced
        Euler characteristic of some simplicial complex on \(n\) vertices.
        Bj\"orner and Kalai
        \cite{extEP} showed that the absolute value of the reduced
        Euler characteristic of a simplicial complex on \(n\) vertices is 
        \(\le \binom{n-1}{\lfloor \frac{n-1}{2} \rfloor}\).
      \end{proof}

      \begin{corr}\label{corr:varmax}
        Let \(f =\sum_{m \in [X]} c_m m \in \Z[[X]]\) be the \([X]\)-graded
        Hilbert numerator of a monoid ideal. Let \(m = x_1^{\alpha_1} \cdots
        x_n^{\alpha_n}\), and let \(r\) be the number of \(1 \le i \le n\)
        such that \(\alpha_i > 0\). Then 
        \begin{equation}
          \label{eq:reuc}
          \mathrm{abs}(c_m) \le \binom{r-1}{\lfloor \frac{r-1}{2} \rfloor}
        \end{equation}
      \end{corr}
      \begin{proof}
        Let \(\sigma\) be a permutation of \(X\). Define
        \(\sigma(x_1^{a_1}\cdots x_\ell^{a_\ell}) = \prod_{i=1}^\ell
        x_{\sigma(i)}^{a_i}\), and \(\sigma(\sum_{m \in [X]} c_m m) = \sum_{m
          \in [X]} c_m \sigma(m)\). We let \(\sigma\) act on monoid ideals in
        \([X]\) in the obvious way.
        Then \(\mu\) and \(\nu\) are fix-points
        for the action of \(\sigma\) on \(\Z[[X]]\), and
        \(\sigma(\Characteristic(I)) = \Characteristic(\sigma(I))\) for all
        monoid ideals \(I\). Hence 
        \begin{align*}
          p(\sigma(I)) &= \mu (\nu - \Characteristic(\sigma(I))) \\
          &=  \mu (\sigma \nu - \sigma (\Characteristic(I))) \\ 
          &= \sigma(\mu (\nu - \Characteristic(I))) \\
          &=\sigma (p(I)).
        \end{align*}
        Let \(i_1,\dots,i_r\) be the support of \(m\), that is, \(\alpha_{i_1} 
        >0, \dots \alpha_{i_r} > 0\), and let \(\sigma\) be a permutation
        which takes \(i_1\) to \(1\), \(i_2\) to \(2\), and so on. Then
        \(\sigma(m) = x_1^{\alpha_{\sigma^{-1}(1)}} \dots
        x_r^{\alpha_{\sigma^{-1}(r)}}\), and 
        \begin{displaymath}
          c_m = [m] f = [\sigma(m)]\sigma(f),
        \end{displaymath}
        hence the result follows by applying Proposition~\ref{prop:varmax}.
      \end{proof}

    \end{subsubsection}

  \end{subsection}
\end{section}

\begin{section}{The subring \protect\(\protect\mathcal{S}\protect\),
    \protect\lfg\ ideals, and their generalised Hilbert numerators}
  For this section, we fix a positive integer \(r\) and set-partition \(Y\) of 
  the set of variables:
  \(X= \cup_{\ell = 1}^r Y_\ell\). There is an
  associated map \(y: \Nat^+ \to \set{1,\dots,r}\) such that \(x_n \in
  Y_{y(n)}\).
  We denote by
  \(\deg\) the associated \(r\)-multi-grading, that is, the monoid
  homomorphism 
  \begin{equation}
    \label{eq:deg}
    \begin{split}
    \deg: [X] & \to \Nat^r \\
    x_i & \mapsto e_{y(i)}\\
    x_1^{\alpha_1} \cdots x_n^{\alpha_n} & \mapsto \alpha_1 \deg(x_1) + \cdots
    \alpha_n \deg(x_n)
    \end{split}
  \end{equation}
  where \(e_1, \dots, e_r\) are the natural basis elements of \(\Nat^r\).
  In particular, if \(r=1\), then \(\deg(m) = \tdeg{m}\).
  Note that, since \(r\) is finite, \RP\ is indeed \(\Nat^r\)-graded by means
  of \(\deg\), even if
  it is not \([X]\)-graded. We say that an ideal is \(r\)-homogeneous if it is 
  homogeneous with respect to this grading. Clearly, \(r\)-homogeneous ideals
  are homogeneous, and all monomial ideals are \(r\)-homogeneous. Furthermore
  we have: 
  
  \begin{proposition}
    Let \(I\) be an ideal of \RP. Then the following are equivalent:
    \begin{enumerate}[(i)]
    \item \(I\) is \(r\)-homogeneous and \alfg,
    \item \(I\) can be generated by \(F= \cup_{\alpha \in \Nat^r} F_\alpha\),
      where each \(F_\alpha\) is a finite set of \(r\)-homogeneous elements of 
      multi-degree \(\alpha\).
    \item For each \(\alpha \in \Nat^r\), 
      \begin{equation}
        \label{eq:kspace}
        \dim_K \left(
        \frac{I_\alpha}{\sum_{\substack{\beta + \gamma = \alpha \\ \beta,
              \gamma \neq 0}} \RP_\beta I_\gamma}
        \right) < \infty
      \end{equation}
    \end{enumerate}
  \end{proposition}
  \begin{proof}
    For each total degree \(d\), there are only finitely many multi-degrees in 
    \(\Nat^r\) of total degree \(d\). Thus \((i)\) and \((i)\) are equivalent.
    The equivalence of \((ii)\) and \((iii)\) is parallel to
    Theorem~\ref{thm:lfg} and is proved in the same way (see
    \cite{Snellman:GbRprimPubl}). 
  \end{proof}
  
  \begin{definition}
    Denote by \(\mathcal{S} \subset \Z[[X]]\) the subring consisting of all
    \(f \in \Z[[X]]\) fulfilling the equivalent conditions below:
    \begin{enumerate}
    \item \(f = f_0 + f_1 + f_2 + f_3 +\cdots\) with \(f_i \in \Z[X]_i\),
    \item \(f(t,t,t,\dots)\), the substitution of each \(x_i\) with the new
      formal indeterminate \(t\), is defined,
    \item \(f = \sum_{\alpha \in \Nat^r} f_\alpha\) with \(f_\alpha \in
      \Z[X]_\alpha\), 
    \item \(f(t_{y(1)},t_{y(2)},t_{y(3)},\dots)\), the substitution of each
      \(x_i\) with the new 
      formal indeterminate \(t_{y(i)}\), is defined,
    \end{enumerate}
    Denote the map \(\mathcal{S} \ni f \mapsto
    f(t_{y(1)},t_{y(2)},t_{y(3)},\dots) \in \Z[[t_1,\dots,t_r]]\) by
    \(\Collapse = \Collapse^y\). 
  \end{definition}

  \begin{theorem}\label{thm:plfg}
    Let \(I \subset [X]\) be a monoid ideal. Then the following are
    equivalent:
    \begin{enumerate}
    \item \(p(I) \in \mathcal{S}\),
    \item \(I\) is \alfg,
    \item \(w(I) \in \mathcal{S}\).
    \end{enumerate}
  \end{theorem}
  \begin{proof}
    Write
    \begin{align*}
      w(I) &= 0 + w_1 + w_2 + w_3 + \cdots, \qquad \tdeg{w_i} = i \\ p(I) &= 1
      + p_1 + p_2 + p_3 + \cdots, \qquad \tdeg{p_i} = i.
    \end{align*}
    It is immediate that \(I\) is \alfg\ if and only if \(w(I) \in
    \mathcal{S}\), which occurs precisely when every \(w_i\) is a polynomial.
    We note that if \(\sigma \subset W(I)\) has cardinality \(u\), and the
    minimal and maximal total degree of elements in \(\sigma\) is \(c\) and
    \(d\), respectively, then \(c + 1 \le \tdeg{\lcm(\sigma)} \le
    ud\). Clearly, the terms of \(w(I)\) contributing to \(p_i\) in
    \eqref{eqn:sigma} have total degree \(\le i\).
    
    Hence, if \(I\) is \alfg, so that each \(w_i\) is a polynomial, then only
    the various lcm's of elements in the support of \(w_1, \dots, w_d\) may
    contribute to \(p_d\). The number of elements in the support of \(w_d\) is
    thus \(\le 2^{(\# w_1 + \cdots \#w_d)} < \infty\).
    
    Conversely, if \(I\) is not \alfg, suppose that \(w_1, \dots, w_d\) are
    polynomials, but that \(w_{d+1}\) is not. Using \eqref{eqn:sigma} we see
    that \(p_{d+1}\) receives contribution from a finite number of terms
    stemming from lcm's of elements in the support of \(w_1, \dots, w_d\), and
    from the non-polynomial \(w_{d+1}\). Thus \(p_{d+1}\) is not a polynomial.
  \end{proof}

  \begin{corr}\label{corr:lpcond}
    Let \(f =\sum_{m \in [X]} c_m m \in \Z[[X]]\). Then \(f \in
    p(\ideal{lm})\) 
    if and only if the following conditions hold:
    \begin{enumerate}
    \item \(\forall m \in [X]: \quad \sum_{\divides{s}{m}} c_s \in
      \set{0,1}\) ,
    \item If \(\sum_{\divides{s}{m}} c_s =1\) and \(\divides{t}{m}\) then
      \(\sum_{\divides{s}{t}} c_s = 1\),
    \item \(f \in \mathcal{S}\).
    \end{enumerate}
  \end{corr}
  \begin{proof}
    This follows from Proposition~\ref{prop:pcond} and
    Theorem~\ref{thm:plfg}. 
  \end{proof}

  We henceforth regard \(\mathcal{S}\) as a topological ring having the
  topology given by the total degree filtration. We have that this topology is 
  the same as the one given by any \(r\)-multi degree filtration
  in the sense that if \(f_n \to f\) if for each multi-degree \(\alpha\),
  there is an \(N(\alpha)\) so that \(f_n\) and \(f\) agrees in multi-degree
  \(\le \alpha\) whenever \(n \ge N(\alpha)\): here \(\le \alpha\) is with
  respect to some term-order on \(\Nat^r\) which refines the total-degree
  partial order. Similarly, the \(\hat{d}\)-metric on homogeneous ideals gives
  the same topology as an analogous \(r\)-multigraded metric.
  
  \begin{lemma}\label{lemma:Sclos}
    \(\mathcal{S}\) is a closed subset of \(\Z[[X]]\). 
  \end{lemma}
  \begin{proof}
    Suppose that \(f_i \to f\), where \(f_i \in \mathcal{S}\). Fix a total
    degree \(d\). By the definition of the total degree filtration topology,
    there exists an \(N\) such that \(\tdegfilt{d}f_i = \tdegfilt{d}f\) for
    all \(i > N\). Since for all \(f_i\), \(\tdegfilt{d}f_i\) is a polynomial, 
    this is true for \(\tdegfilt{f}\), as well.
  \end{proof}

  \begin{theorem}\label{thm:mScont}
    \(\Collapse: \mathcal{S} \to \Z[[t_1,\dots,t_r]]\) is continuous and
    clopen, when 
    \(Z[[t_1,\dots,t_r]]\) is 
    given the \((t_1,\dots,t_r)\)-adic topology.
  \end{theorem}
  \begin{proof}
    We assume for simplicity that \(r=1\).
    Suppose that \(f_i \to f\) in \(\mathcal{S}\). Fix an integer \(d\), and
    choose an \(N(d)\) such that \(f_i -f \in \bar{\mathfrak{m}^{d}}\) for \(i
    > N(d)\). Thus for \(i > N(d)\) we have that the \(t^d\) coefficient of
    \(\Collapse(f_i)\) and \(\Collapse(f)\) coincides. This shows that
    \(\Collapse(f_i) \to \Collapse(f)\).

    To show that this map is clopen, we pick a basic clopen subset 
    \(O_{f,d} = 
    \setsuchas{g \in \mathcal{S}}{\tdegfilt{d} f = \tdegfilt{d} g}\), where
    \(d\) is a positive integer, and \(f \in \mathcal{S}\). 
    Then \(\Collapse(O_{f,d}) = \setsuchas{h \in \Z[[t]]}{\tdegfilt{d} h =
      \tdegfilt{d} \Collapse(f)}\), and this is a basic clopen set of
    \(\Z[[t]]\). 
  \end{proof}

  \begin{lemma}\label{lemma:charhom}
    The characteristic function is a continuous mapping from the set of \alfg\
    monomial ideals, with the \(\hat{d}\) metric, to \(\mathcal{S}\). In fact,
    it is a homeomorphism onto its image.
  \end{lemma}
  \begin{proof}
    Obvious.
  \end{proof}

\begin{lemma}\label{lemma:charclosed}
    The set of characteristic functions of \alfg\ monoid ideals is a closed subset of
    \(\mathcal{S}\) (and of \(\Z[[X]]\)).
  \end{lemma}
  \begin{proof}
    This follows from the previous Lemma and from Lemma~\ref{lemma:Sclos}.
  \end{proof}

  \begin{lemma}\label{lemma:pclosed}
    The set \(p(\ideal{lm}) \subset \mathcal{S}\) is closed.
  \end{lemma}
  \begin{proof}
    By the previous Lemma, the set of characteristic functions of \alfg\ monoid
    ideals is a closed subset of  \(\Z[[X]]\). By Lemma~\ref{lemma:transclos}, 
    the mapping \(f \mapsto \mu(\nu - f)\) is a closed mapping, hence
    \(p(\ideal{lm})\) is a closed subset of \(\Z[[X]]\). From
    Theorem~\ref{thm:plfg} we have that \(p(\ideal{lm}) \subset
      \mathcal{S}\), hence it is closed in there.
  \end{proof}

  We henceforth assume that \(\ideal{lm}\) have the \(\hat{d}\)-topology.

  \begin{proposition}\label{prop:contcon}
    Let \(I,I_1,I_2,I_3,\dots\) be \alfg\ monomial ideals in \RP. The following
    are 
    equivalent:
    \begin{enumerate}
    \item \(\hat{d}(I_n,I) \to 0\),
    \item \(\Characteristic(I_n) \to \Characteristic(I)\),
    \item \(q(I_n) \to q(I)\),
    \item \(p(I_n) \to p(I)\),
    \item \(w(I_n) \to w(I)\).
    \end{enumerate}
    Furthermore, if the conditions above are satisfied, then
    \begin{displaymath}
      \Collapse(p(I_n)) \to \Collapse(p(I)).
    \end{displaymath}
  \end{proposition}
  \begin{proof}
    By the previous lemma, \(I_n \to I\) if and only if \(\Characteristic(I_n)
    \to \Characteristic(I)\). Since the endomorphism given by multiplication
    with a fixed element in a topological ring is continuous,
    \begin{displaymath}
      \Characteristic(I_n) \to \Characteristic(I) \quad \iff \quad q(I_n) \to
      q(I) \quad \iff \quad p(I_n) \to p(I).
    \end{displaymath}
    If \(w(I_n) \to w(I)\), then fixing a total degree \(d\), we get that
    there exists an \(N(d)\) such that \(w(I_n) - w(I) \in
    \hat{\mathfrak{m}^d}\) for \(n \ge N(d)\).  It then follows that
    \(\Characteristic(I_n) - \Characteristic(I) \in \hat{\mathfrak{m}^d}\) for
    \(n \ge N(d)\).  The converse also holds.
    
    The last assertion follows immediately from the fact that \(\Collapse\) is
    continuous.
  \end{proof}

  We now recall a theorem by Macaulay, which says that if \(I \subset K[X_n]\)
  is a homogeneous ideal, and \(>\) is a term-order on \([X_n]\), then
  \(\frac{K[X_n]}{I}\) and \(\frac{K[X_n]}{\init_>(I)}\) have the same
  \(\Nat\)-graded Hilbert series  (see for instance \cite{Ebud:View}). It is
  also true that if \(I\) is \(r\)-multigraded, when \(K[X_n]\) is
  \(\Nat^r\)-graded using the partition \(Y \cap X_n\), then the above
  algebras have in fact the same \(\Nat^r\)-graded Hilbert series.
  Using this, and our previous
  results, we get:

  \begin{theorem}
    Suppose that \(>\) is a term-order on \([X]\).  Let \(I \subset \RP\) be a
    \(r\)-homogeneous \alfg\ ideal, and define \(g_n \in \Z[t_1,\dots,t_r]\) by
    requiring that 
    \begin{displaymath}
      \frac{g_n}{\prod_{i=1}^n (1-t_{y(i)})}
    \end{displaymath}
    is the \(\Nat^r\)-graded Hilbert series of
    \begin{displaymath}
      \frac{K[X_n]}{\restrict[n]{I}}.
    \end{displaymath}
    Then, \(g_n \to \Collapse(p(I))\) as \(n \to \infty\), and
    \(\Collapse(p(I)) \in \Z[[t_1,\dots,t_r]]\).
  \end{theorem}
  \begin{proof}
    From Theorem~\ref{thm:inconv} we know that
    \[\hat{d}(\init_>(\restrict[n]{I})^e, \init_>(I)) \to 0.\]
    Then Proposition~\ref{prop:contcon} gives that
    \begin{displaymath}
      p(\init_>(\restrict[n]{I})^e) \to p(\init_>(I)),
    \end{displaymath}
    hence, using Theorem~\ref{thm:mScont}, we have that
    \begin{displaymath}
      \Collapse(p(\init_>(\restrict[n]{I})^e)) \to \Collapse(p(\init_>(I))).
    \end{displaymath}
    It is clear that
    \begin{displaymath}
      p^n(\init_>(\restrict[n]{I})) = p(\init_>(\restrict[n]{I})^e).
    \end{displaymath}
    As we remarked above, a \(r\)-homogeneous ideal have the same
    \(\Nat^r\)-graded 
    Hilbert series as its initial ideal, so
    \begin{displaymath}
      g_n = \Collapse(p^n(\init_>(\restrict[n]{I}))),
    \end{displaymath}
    hence
    \begin{displaymath}
      g_n \to \Collapse(p(\init_>(I))).
    \end{displaymath}
  \end{proof}

  In \cite{Snellman:HilbNumPubl}, the result above (for \(r=1\)) was proved
  through a 
  different route, and the power series \(\Collapse(p(I))\) was called the
  \emph{generalised Hilbert numerator} of \(I\).
  
  We note two simple corollaries:
  \begin{corr}
    If \(I,J\) are \(r\)-homogeneous \alfg\ ideals of \RP, and if
    \(\restrict[n]{I}\) and 
    \(\restrict[n]{J}\) have the same \(\Nat^r\)-graded Hilbert series for all
    \(n\), then \(I\) 
    and \(J\) have the same \(r\)-graded Hilbert numerator.
  \end{corr}
  
  \begin{corr}\label{corr:monsuf}
    If \(I\) is an \(r\)-homogeneous \alfg\ ideal of \RP, and if \(>\) is a
    term-order on \([X]\), 
    then \(I\) and \(\init_>(I)\) have the same \(r\)-graded Hilbert numerator.
  \end{corr}
  In particular, \(\Collapse(p(\ideal{l})) = \Collapse(p(\ideal{lm}))\), that
  is, all \(r\)-graded Hilbert numerators of \alfg\ ideals can be obtained from 
  \alfg\ monomial ideals.
  
  \begin{subsection}{Polynomial \protect\(r\protect\)-graded Hilbert
      numerators} 
    \begin{definition}
      We put \(\mathcal{T}^Y = \Collapse^{-1}(\Z[t_1,\dots,t_r])\).
      If \(p(I) \in \mathcal{T}^Y\) we say that \(I\) has \emph{polynomial
      \(r\)-graded Hilbert numerator}.
    \end{definition}
    
    \begin{lemma}
      Suppose that \(Y'\) refines \(Y\). Denote by \(0\) the partition
      \(X=X\). Then 
      \begin{math}
        \Z[X] \subsetneq  \mathcal{T}^{Y'} \subset
        \mathcal{T}^Y \subset \mathcal{T}^0 \subsetneq \mathcal{S}. 
      \end{math}
    \end{lemma}
    \begin{proof}
      The inclusions are obvious. To see that the strict ones are indeed
      strict, consider the
      following examples:
      \begin{math}
        \sum_{i=1}^\infty (x_1^i - x_2^i) \in \mathcal{T}^0 \setminus \Z[X],
      \end{math}
      \begin{math}
        \sum_{i=1}^\infty x_1^i \in \mathcal{S} \setminus \mathcal{T}^Y.
      \end{math}
    \end{proof}
    
    \begin{example}\label{ex:23gen}
      There are \alfg\ monomial ideals which have Hilbert numerators in
      \(\mathcal{S} \setminus \mathcal{T}^0\).
      Let \(I\) be generated by \(a_i = x_1x_2x_3 \cdots x_{i-1}x_i^2\), for
      \(i \ge 1\), and \(b_j = x_1 x_2 x_3 \cdots x_{j-2} x_j^6\), for
      \(j \ge 2\). Put \(p_n = p^n(\restrict[n]{I})\) for \(n > 0\), and
      define \(p_0  = 1\). We claim that 
      \begin{equation}
        \label{eq:23rec}
        p_n = p_{n-1} + (-1)^n v_n, 
        \qquad v_n = (x_{n-1} + x_n^4)x_1x_2\cdots
        x_{n-2}x_n^2 \prod_{i=1}^{n-1}(x_i - 1)
      \end{equation}
      To see this, we first note that 
      \begin{math}
        (x_{n-1} + x_n^4)x_1x_2\cdots x_{n-2}x_n^2 = a_n + b_n.
      \end{math}
      By induction, we have that
      \begin{math}
        (x_{n-1} + x_n^4)x_1x_2\cdots x_{n-2}x_n^2 
        \prod_{i=1}^{n-1}(x_i - 1) 
      \end{math}
      consists of those monomials which can be formed as a lcm of 
      \(\set{a_n} \cup S\) or of \(\set{b_n} \cup S\) or of 
      \(\set{a_n,b_n} \cup S\), where 
      \(S \subset  \set{a_1,a_2,\dots,a_{n-1},b_2,b_3,\dots,b_{n-1}}\). 
      Note that monomials
      \(x_1^{\alpha_1} \cdots x_n^{\alpha_n}\) 
      with \(\alpha_i = \alpha_j =  6\) for \(i < j\) does not occur. 
      This can be readily explained:
      every such monomial can be expressed as a lcm in two different ways, by
      either including or omitting the superfluous generator \(a_i\).
      
      Hence, it follows that 
      \begin{equation}
        \label{eq:23exp}
        \begin{split}
        p_n &= 1 + \sum_{i=1}^n (-1)^i v_i \\
        \lim p_n = p & = 1 + \sum_{i=1}^\infty (-1)^i v_i
        \end{split}
      \end{equation}
      We have that \(p \in \mathcal{S}\setminus
      \Z[X]\). As we shall see, \(p \in \mathcal{T}\).

      Setting each \(x_i=t\) in \eqref{eq:23exp} we have that
      \begin{equation}
        \label{eq:23t}
        \begin{split}
          p(t) &= 
          1 - t^2 + \sum_{n=2}^\infty (-1)^n (t-1)^{n-1}(t+t^4)t^n \\
          &=  1 - t^2 + \frac{t+t^4}{t-1} \sum_{n=2}^\infty (-1)^n (t-1)^n
          t^n\\ 
          & = 1 - t^2 + \frac{t+t^4}{t-1} \left( \frac{1}{1+t(t-1)} - 1 +
            (t-1)t\right) \\ 
          & = 1 - t^2 - t^3 + t^5
        \end{split}
      \end{equation}
    \end{example}
    
    \begin{lemma}\label{lemma:Q}
      Let \(r=1\).
      Let \(d,a_1,\dots,a_d\) be integers, with \(d > 0\). Then the set
      of all \(\tdeg{p(I)}\),
      where \(I\)  is finitely generated and generated
          in degrees \(\le d\), and has 
          \begin{equation}
            \label{eq:q}
          \Collapse(p(I)) = 1 + a_1+ \dots
          a_dt^d + O(t^{d+1}), 
          \end{equation}
        is either empty, or has a maximum, which we denote by
        \(Q_d(a_1,\dots,a_d)\). 
    \end{lemma}
    \begin{proof}
      We claim that there are positive integers \(u_1,\dots,u_d\) such that if 
      \(I\) is a finitely generated monomial ideal generated in degrees \(\le
      d\) satisfying \eqref{eq:q}, then \(\Collapse(w(I)) = w_1t + \dots
      w_dt^d\) with \(w_i \le u_i\) for \(1 \le i \le d\). Assuming the claim,
      it is clear that the  total degree of \(p(I)\) is \(\le \sum_{i=1}^d i
      u_i\), since this is a bound of the lcm of all the generators.

      To establish the claim, we note that 
      \(a_1 = - w_1\), and assume by
      induction that we have shown that \(u_1,\dots,u_i\) exist. 
      We note that the minimal generators which affect \(a_{i+1}\) are those
      of degree 
      \(i+1\), which each contribute with \(-1\), and also \(s\)-tuples
      \(m_1,\dots,m_s\) with \(m_\ell \in W(I)\), \(\tdeg{m_\ell} < i+1\),
      and with 
      \(\tdeg{\lcm(m_1,\dots,m_s)} = i+1\), which each contribute \((-1)^s\).
      If we pick \(\lambda_1\) elements of \(W(I)_1\), \(\lambda_2\) elements
      of \(W(I)_2\), et cetera, then for the resulting lcm to be of total
      degree \(i+1\) it is necessary that \(\lambda_1 + \lambda_2+ \dots +
      \lambda_i\ge 
      i+1\) and that \(\lambda_\ell < i+1\) for all \(\ell\); thus only
      finitely many \(\lambda = (\lambda_1,\dots,\lambda_i)\) are relevant.

      We thus have that
      \begin{equation}
        \label{eq:range}
        a_{i+1} =  -w_{i+1} + \sum_{\lambda} (-1)^{c(\lambda)}R_\lambda,
      \end{equation}
      where \(\lambda=(\lambda_1,\dots,\lambda_i)\),
      \(\tdeg{\lambda} \ge i+1\), \(0 \le \lambda_\ell < i+1\) for \(1 \le
      \ell  
      \le j \le i\),   \(c(\lambda) \) is the number of non-zero entries in
      \(\lambda\). 
      The symbol \(R_\lambda\) denotes a finite interval \([0,L]\) of
      integers,  
      where \(L\) is the maximal numbers of lcm's of \(\lambda_1\) elements of
      degree 1, drawn from a set of cardinality \(u_1\), \(\lambda_2\)
      elements of degree 2, drawn from a set of cardinality \(u_2\), and so
      on,  which have total 
      degree \(i+1\).
      For instance, if \(i=1\) and \(\lambda = (2)\) then
      \(L=\binom{u_1}{2}\), if \(i=2\) and \(\lambda =(1,1)\) then \(L=u_1
      u_2\).  

      These finite intervals are added using interval
      arithmetic, so that \([a,b] + [c,d] = [a+c,b+d]\). We can the deduce
      that
      \begin{equation}
        \label{eq:range2}
        w_{i+1} = - a_{i+1} + \sum_{\lambda} (-1)^{c(\lambda)}R_\lambda =
        [-C,D] 
      \end{equation}
      for some integers \(C,D\). Putting \(u_{i+1} = D\) we have the desired
      bound. 
    \end{proof}
    
    \begin{lemma}\label{lemma:truncHilb}
      Let \(r=1\) and suppose that \(f(t) \in \Collapse(p(\ideal{l}))\), in
      other words, 
      that \(f(t)\) is 
      the \(\Nat\)-graded Hilbert numerator of some \alfg\ ideal.
      Suppose that \(f(t) = 1 + a_1t +
      \cdots + a_dt^d + t^{d+r+1}g(t)\), with \(r > Q_d(a_1,\dots,a_d)\).
      Then 
      \(1 + a_1t + \cdots + a_dt^d \in \Collapse(p(\ideal{l}))\).
    \end{lemma}
    \begin{proof}
      Let \(f = \Collapse(p(I))\), where \(I\) is a \alfg\ monomial ideal.
      Let \(I_{\le d}\) denote the  ideal generated by everything in \(I\) of
      total  degree \(\le d\).
      Since the
      maximal degree of a lcm of  
      the generators of degree \(\le d\) is \(Q_d(a_1,\dots,a_r)\), it follows
      from \eqref{eqn:sigma} and Lemma~\ref{lemma:Q} that
      \(\Collapse(p(I_{\le d}) = 1 + a_1t + \cdots + a_dt^d\). 
    \end{proof}

    \begin{theorem}\label{thm:poly}
      Let \(r=1\). If \(I \subset [X]\) is a \alfg\ monoid ideal with \(p(I)
      \in \mathcal{T}\), then there exists a positive integer \(N\) and a
      monoid ideal \(J \subset [X_N]\) so that \(\Collapse(p(J^e)) =
      \Collapse(p(I))\). 
    \end{theorem}
    \begin{proof}
      Let \(f = \Collapse(p(I)) =1+a_1t +a_2t^2 +  \dots a_dt^d\), 
      and let \(f_n = \Collapse(p^n(\restrict[n]{I}))\). Then \(f_n \to f\) in
      \(\Z[[t]]\), 
      with respect to the \((t)\)-adic topology.  Let \(r >
      Q_d(a_1,\dots,a_d)\), and choose \(N\) such that for
      \(n \ge N\), \(f_n - f \in (t^r)\). Then Lemma~\ref{lemma:truncHilb}
      shows that there is a monoid ideal \(J\) in \([X_n]\) with \(f\) as its
      \(\Nat\)-graded Hilbert numerator.
    \end{proof}
    
    \begin{corr}\label{corr:clos}
      The set of polynomial \(\Nat\)-graded Hilbert numerators of \alfg\
      ideals in \(\RP\)  is equal to
      the set of \(\Nat\)-graded Hilbert numerators of homogeneous ideals in
      finitely many variables. This set is dense in the set of all possible
      \(\Nat\)-graded Hilbert numerators of \alfg\ ideals.
    \end{corr}
    \begin{proof}
      From Theorem~\ref{thm:poly} we get that all polynomial \(\Nat\)-graded
      Hilbert numerators can be obtained from ideals generated in finitely
      many variables. To prove the second assertion, we note that if \(I\) is
      \alfg, \(d>0\), and \(I_{\le d}\) is the ideal generated by everything in 
      \(I\) of degree \(\le d\), then \(\Collapse(p(I)) \simeq
      \Collapse(p(I_{\le d})) \mod (t^{d+1})\), and since \(I_{\le d}\) is
      finitely generated, \(p(I) \in \Z[X]\) hence \(\Collapse(p(I)) \in
      \Z[t]\). 
    \end{proof}

    \begin{theorem}\label{thm:class}
      Let, for every pair of integers \(0 < a \le b\), \(G_{a,b}\) denote the
      set 
      \begin{displaymath}
        \setsuchas{(1-t)^b (1 + a_1 t + a_2 t^2 + a_3 t^3 + \cdots)}{a_1 = a,
          \quad \forall i: 0 \le a_{i+1} \le a_i^{<i>}},
      \end{displaymath}
      where 
      \[u^{<d>} = \binom{k(d)+1}{d+1} + \binom{k(d-1)+1}{d} + \cdots +
      \binom{k(1) +1}{2}\] when \(u\) has \(d\)-th Macaulay expansion
      \[u = \binom{k(d)}{d} + \binom{k(d-1)}{d-1} 
      + \cdots + \binom{k(1)}{1}\] 
      (see \cite{BrunsHerzog:CM}). Then the set of polynomial \(\Nat\)-graded
      Hilbert numerators is \(\cup_{0 < a \le b} G_{a,b}\), and the closure of 
      this set in \(\Z[[t]]\) is exactly the set of \(\Nat\)-graded Hilbert
      numerators of \alfg\ ideals in \RP.
    \end{theorem}
    \begin{proof}
      It follows from a well-know classification by Macaulay (see
      \cite[Theorem 4.2.10]{BrunsHerzog:CM}) that the set of (generating
      functions of) Hilbert functions 
      of homogeneous quotients of polynomial rings with finitely many
      indeterminates is 
      \begin{displaymath}
        \setsuchas{1 + a_1 t + a_2 t^2 + a_3 t^3 + \cdots}{\forall i: 0 \le
          a_{i+1} \le a_i^{<i>}}.
      \end{displaymath}
      The function \(1 + a_1 t + a_2 t^2 + a_3 t^3 + \cdots\) can be realised
      as the Hilbert function of a quotient of \(K[x_1,\dots,x_{a_1}]\) with a 
      monomial ideal; we are of course free to use more variables, if we so
      desire. The first part of the theorem is therefore demonstrated.

      To prove the second part, we proceeds as follows. 
      We know by Lemma~\ref{lemma:pclosed} that 
      \[p(\ideal{lm})= \setsuchas{\mu(\nu - \Characteristic(I))}{I \in
        \ideal{lm}}\]  
      is a closed subset of \(\mathcal{S}\). We have that \(\Collapse:
      \mathcal{S} \to \Z[[t]]\) is a closed map (Theorem~\ref{thm:mScont}),
      hence 
      \(\Collapse(p(\ideal{lm}))\) is closed in \(\Z[[t]]\).

      Now, Corollary~\ref{corr:clos} shows that the set of polynomial
      \(\Nat\)-graded Hilbert 
      numerators is dense in the set of all \(\Nat\)-graded Hilbert
      numerators; since this latter set is closed, the second part of the
      theorem follows.
    \end{proof}

  \end{subsection}
  
\end{section}


\bibliographystyle{plain}
\bibliography{journals,snellman,articles}
\end{document}